\documentclass[11pt]{article}

\usepackage{etex}
\usepackage[margin={1.2in}]{geometry}
\usepackage[titletoc]{appendix}

\usepackage{rotating}

\usepackage{longtable}

\usepackage[sc]{mathpazo}
\usepackage{courier}

\usepackage[utf8]{inputenc}
\usepackage[T1]{fontenc}

\usepackage[english]{babel}
\usepackage{blindtext}

\usepackage{amsmath}

\usepackage{microtype}

\usepackage{amsmath}
\usepackage{amssymb}
\usepackage{amsthm}
\usepackage{color}
\usepackage{latexsym,extarrows}

\usepackage[all]{xy}

\usepackage[stable,multiple]{footmisc}

\RequirePackage[font=small,format=plain,labelfont=bf,textfont=it]{caption}
\begin{document}
\def\e#1\e{\begin{equation}#1\end{equation}}
\def\ea#1\ea{\begin{align}#1\end{align}}
\def\eq#1{{\rm(\ref{#1})}}
\theoremstyle{plain}
\newtheorem{thm}{Theorem}[section]
\newtheorem{lem}[thm]{Lemma}
\newtheorem{prop}[thm]{Proposition}
\newtheorem{cor}[thm]{Corollary}
\theoremstyle{definition}
\newtheorem{dfn}[thm]{Definition}
\newtheorem{ex}[thm]{Example}
\newtheorem{rem}[thm]{Remark}
\newtheorem{conjecture}{Conjecture}

\newtheorem*{mainthmA}{{\bf Theorem A}}
\newtheorem*{maincorC}{{\bf Corollary C}}
\newtheorem*{mainthmB}{{\bf Theorem B}}
\newtheorem*{mainthmD}{{\bf Theorem D}}

\newcommand{\comment}[1]{\textcolor{red}{[#1]}} 

\newcommand{\LD}{\langle}
\newcommand{\RD}{\rangle}

\numberwithin{equation}{section}

\title{\bf  Mirror Symmetry for Plane Cubics Revisited}
\author{Jie Zhou}
\date{}
\maketitle

\begin{abstract}
In this expository note we discuss some arithmetic aspects of the mirror symmetry for plane cubic curves. We also explain how the Picard-Fuchs equation can be used to reveal part of these arithmetic properties. The application of Picard-Fuchs equations in studying the genus zero Gromov-Witten invariants of more general Calabi-Yau varieties and the Weil-Petersson geometry on their moduli spaces will also be discussed.
\end{abstract}

\setcounter{tocdepth}{2} \tableofcontents

\section{Introduction}
\label{secintro}

Mirror symmetry is a surprising and yet to be further explored symmetry on the \emph{moduli space} of Calabi-Yau (CY) varieties. 
It is usually referred to in the form of the slogan "Mirror symmetry exchanges the \emph{K\"ahler structure} of a
Calabi-Yau space with the \emph{complex structure} of its mirror."\\

The simplest CY variety is the elliptic curve for which the above mirror phenomenon is mostly understood, see \cite{Dijkgraaf:1995}.
As a complex variety, an elliptic curve has the form $E_{\tau}=\mathbb{C}/(\mathbb{Z}\oplus \mathbb{Z}\tau)$ with $\tau\in \mathcal{H}$, where $\mathcal{H}$ is the upper-half plane.
The parameter $\tau$ determines the "shape" of the lattice $\Lambda_{\tau}=\mathbb{Z}\oplus \mathbb{Z}\tau$ and captures the complex structure. 
A K\"ahler structure is described  by
\begin{equation}
\omega=r\omega_{*}\,,\quad \omega_{*}= {\sqrt{-1}\over 2} {1\over \mathrm{ Im} \tau}dz_{\tau}\wedge d\overline{z}_{\tau} \,.
\end{equation}
Here $z_{\tau}$ is the standard complex coordinate system on (the universal cover of) $E_{\tau}$\,.
The K\"ahler form $\omega_{*}$, satisfying $\int_{E_{\tau}}\omega_{*}=1$, is taken to be the basis for the tangent space of the space of K\"ahler structures.
The parameter $r$ is a real positive number.

To make the space of K\"ahler structures a complex variety, one 
introduces the "B-field" $B\in H^{1,1}(E_{\tau},\mathbb{C})$
and considers the space of "complexified" K\"ahler structures of the form 
\begin{equation}
\omega_{\mathbb{C}}=B+\sqrt{-1}\omega:=t\omega_{*}\,,
\quad t=\theta+\sqrt{-1} r\,.
\end{equation}
The condition $r> 0$ gets translated into the condition $\mathrm{Im}\,t>0$ on the "complexified size" $t$. 

Now we can describe the CY structure of an elliptic curve in terms of the two parameters $(t,\tau)\in \mathcal{H}\times \mathcal{H}$. 
We denote the corresponding CY structure on the elliptic curve by $E_{t,\tau}$.
Hereafter, to avoid potential confusion, we shall use
$\mathcal{H}_{\mathrm{K}}$ to denote the first copy of $\mathcal{H}$
and $\mathcal{H}_{\mathrm{C}}$ the second one.\\

Mirror symmetry says that the mirror $\check{E}_{\check{t},\check{\tau}}$ of $E_{t,\tau}$ is given by $E_{\tau,t}$.
Therefore, it is a tautology that the space of complexified K\"ahler structures of $E$
is identified with the space of complex structures of $\check{E}$ and vice versa. It is in this sense that mirror symmetry exchanges the K\"ahler structure of a CY variety with the complex structure of its mirror. 

Among many other things, it also conjectures that the K\"ahler geometry (A-model)
of an elliptic curve
should be equivalent to the complex geometry (B-model) of its mirror elliptic curve, and vice versa.
For example, the computation of genus zero \emph{Gromov-Witten invariants}
in the A-model can be translated into the study of \emph{variation of Hodge structures} in the B-model \cite{Candelas:1990rm}.
The same story is conjectured to be true for general CY varieties but a complete and conceptual understanding is still lacking, see \cite{Cox:2000vi, Hori:2003ic} for a review on this subject.\\

The above two spaces $\mathcal{H}_{\mathrm{K}}$ and $\mathcal{H}_{\mathrm{C}}$ are actually not the desired "moduli" spaces since
a lot of points in each space should be identified. 
The space $\mathcal{H}_{\mathrm{C}}$ consisting of the $\tau$'s can be naturally regarded as the moduli space of complex structures with \emph{markings}, or alternatively the moduli space of marked polarized Hodge structures.
There is a natural $\mathrm{SL}_{2}(\mathbb{Z})$-action on $\mathcal{H}_{\mathrm{C}}$ which identifies different markings, yielding the quotient $\mathcal{M}_{\mathrm{C}}=\mathrm{SL}_{2}(\mathbb{Z})\backslash \mathcal{H}_{\mathrm{C}}$
as the true moduli space\footnote{This is not a fine moduli space due to the existence of torsion elements in the group but only a coarse moduli space. Strictly speaking this space should be treated as an orbifold or even more generally a stack.}of complex structures of the elliptic curve.
By going from $\mathcal{H}_{\mathrm{C}}$ to $\mathcal{M}_{\mathrm{C}}$ one simply forgets about the markings.\\

To meet the expectation from mirror symmetry, there should exist an $\mathrm{SL}_{2}(\mathbb{Z})$-action on $\mathcal{H}_{\mathrm{K}}$ as the mirror of the $\mathrm{SL}_{2}(\mathbb{Z})$-action on $\mathcal{H}_{\mathrm{C}}$.
While one of the generators of $\mathrm{SL}_{2}(\mathbb{Z})$ given by $T:t\mapsto t+1$ can be easily realized by thinking of the B-field as valued in 
$H^{2}(E_{\tau},\mathbb{C})/H^{2}(E_{\tau},\mathbb{Z})$, it is not easy to find
a direct geometric interpretation of the other generator $S:t \mapsto -1/t$.
\footnote{In physics, this follows from the T-duality on the worldsheet.}\\

Another motivation of studying these actions comes from understanding the quasi-modularity in the Gromov-Witten theory
of the elliptic curve.
It has been known that \cite{Dijkgraaf:1995, Kaneko:1995, Bloch:2000, Eskin:2001, Okunkov:2002, Roth:2010ffa} the generating series of simply branched coverings (equivalently Hurwitz numbers)
correspond to the Fourier expansions of some quasi-modular forms \cite{Kaneko:1995} for the modular group $\mathrm{SL}_{2}(\mathbb{Z})$.
In order to describe the modular group action on various constructions (e.g., Hurwitz moduli spaces) in the enumerative geometry which leads to the quasi-modularity, 
it seems necessary to understand the action of the $S$-transform on the parameter $t$ first.
Note that things would become much more clear on the mirror B-model \cite{Dijkgraaf:1995, Bershadsky:1993ta, Bershadsky:1993cx, Aganagic:2006wq, Li:2011thesis, Li:2011b, Li:2011mi, Costello:2012, Bohem:2015} where quasi-modularity is regarded as certain equivariance under the action of $\mathrm{SL}_{2}(\mathbb{Z})$ on the space $\mathcal{H}_{\mathrm{C}}$.

\subsection*{Structure of the note}

This expository note is aimed at finding a conceptual understanding of the $\mathrm{SL}_{2}(\mathbb{Z})$-action on the space of K\"ahler structures of the elliptic curve.
We shall discuss some \emph{arithmetic structures} that are involved in the mirror symmetry of elliptic curves.
These arithmetic structures are studied very little (comparing to the geometric structures) in the literature. We hope that by revealing how they might come into the play for the elliptic curve case
can shed some light on the studies of the mirror symmetry phenomenon in general.

It is believed that in order to have a through understanding of the $\mathrm{SL}_{2}(\mathbb{Z})$-action, the space of complexified K\"ahler structures should be replaced by the space of suitable stability conditions \cite{Bridgeland:2007stability}.
In this note, we shall however only present some mostly speculative discussions in elementary terms.\\

In Section \ref{seclattice} we study the \emph{lattices} that are involved in the mirror symmetry of elliptic curves
which are responsible for the origins of the $\mathrm{SL}_{2}(\mathbb{Z})$-actions. In Section \ref{sectorsion} we discuss the role of the \emph{torsion} in the lattices by studying the mirror of plane cubic curves as an example. 
In Section \ref{secdualities} we explain how the analytic continuation
of the solutions to Picard-Fuchs equations can be used
to detect part of the torsion structure and also reveal some \emph{dualities} between different theories.
Section \ref{secYukawa} is devoted to discussing the genus zero Gromov-Witten invariants of Calabi-Yau varieties and the \emph{Weil-Petersson geometry} on their moduli spaces by making use of the Picard-Fuchs equations.

\subsection*{Acknowledgment} This note is partially based on the talks the author gave at the FRG Workshop "Recent Progress in String Theory and Mirror Symmetry" in May 2015 at Brandeis University, the workshop "Uniformization, Riemann-Hilbert Correspondence, Calabi-Yau Manifolds, and Picard-Fuchs Equations" in July 2015 at  Institut Mittag-Leffler, and the Workshop "Algebraic Varieties" in November 2015 at the Fields Institute. The author would like to thank the organizers for invitations and the institutions for hospitality.  
He is grateful to Kevin Costello, An Huang, Si Li, Lizhen Ji,  Baosen Wu, Shing-Tung Yau and Noriko Yui for their interest and for helpful discussions. He also thanks Yefeng Shen and Siu-Cheong Lau for collaborations on related topics. In addition he thanks the referees for useful comments.

The author is supported by the Perimeter Institute for Theoretical Physics. Research at Perimeter Institute is  supported  by the Government  of Canada  through  Industry Canada and  by the Province of Ontario through the Ministry of Economic Development and Innovation.

\section{Role of (co)homology lattices}
\label{seclattice}

Recall that $\mathcal{H}_{\mathrm{C}}$ is actually the moduli space of complex structures with extra structures on $\check{E}$, namely the markings $\check{m}: \mathbb{Z}\oplus \mathbb{Z}\cong H_{1}(\check{E},\mathbb{Z})$.
Now instead of focusing on how to interpret the space $\mathrm{SL}_{2}(\mathbb{Z}) \backslash\mathcal{H}_{\mathrm{K}}$ as a moduli space, one may ask whether $\mathcal{H}_{\mathrm{K}}$ is the moduli space of K\"ahler structures with certain extra structures on $E$ as well.
A na\"ive attempt is to regard the extra structure as a marking 
\begin{equation}\label{eqndesiredmarking}
m: \mathbb{Z}\oplus \mathbb{Z}\cong \Lambda(E)\,,
\end{equation}
where $\Lambda(E)$ is some rank $2$ lattice constructed from the (co)homology of $E$.

\subsection{Lattice in the A-model}
 
Motivated by the connection to stability conditions and variation of Hodge structures in the B-model, one is immediately led to a natural candidate of the markings in \eqref{eqndesiredmarking} as described below.

We first recall the standard notations from the theory of variation of Hodge structures in the B-model on $\check{E}$.
The Hodge bundle $\mathcal{F}^{0}$ over $\mathcal{H}_{\mathrm{C}}$ is the pull back of the rank two trivial bundle via
$\mathcal{H}_{\mathrm{C}}\rightarrow \mathrm{Gr}(H^{1,0}(\check{E},\mathbb{C}), H^{1}(\check{E},\mathbb{Z})\otimes \mathbb{C)}$
and the Hodge line bundle $\mathcal{F}^{1}$ the pull back of the tautological line bundle.
Locally, the latter has a unique (up to scaling by a holomorphic function) holomorphic section $\Omega$.
A marking $\check{m}$ gives a symplectic basis of cycles $A,B$ in $H_{1}(\check{E}, \mathbb{Z})$. Denote their duals in $H^{1}(\check{E},\mathbb{Z})$, which is the local system that underlies $\mathcal{F}^{0}$, by $\alpha,\beta$ respectively. 
A local trivialization of the Hodge line bundle is 
\begin{equation}
\tau\beta+\alpha\,.
\end{equation}
For a generic section, one has
\begin{equation}\label{eqnBtopform}
\Omega=\pi_{1}\beta+\pi_{0}\alpha\,,
\end{equation}
where $\pi_{0}=\int_{A}\Omega\,, \pi_{1}=\int_{B}\Omega$ are the period integrals with respect to the basis $A,B$.

There is a similar story on the A-model which describes the variation
of quantum cohomology ring, see for example \cite{Hori:2003ic} for a nice account of these discussions. 
The relevant local system is the K-theory group $K_{0}(E)/torsion$ which under the Chern character isomorphism (i.e., tensoring over $\mathbb{Q}$ gives the isomorphism) is identified with $H^{\mathrm{even}}(E,\mathbb{Z})=H^{0}(E,\mathbb{Z})\oplus H^{2}(E,\mathbb{Z})$. 
After tensoring with $\mathbb{C}$, the latter glue to the Hodge bundle $\mathcal{F}^{0}$ over $\mathcal{H}_{\mathrm{K}}$.
The Hodge line bundle $\mathcal{F}^{1}\subseteq \mathcal{F}^{0}$ has a local trivialization
given by
\begin{equation}\label{eqnAlocaltrivialization}
{1\over \sqrt{\mathrm{Td}(E)}}\exp_{q} (\omega_{\mathbb{C}})\,.
\end{equation}
Here the subscript $q$ in $\exp_{q}$ means that in the definition of $\Omega$ as a formal series in $\omega_{\mathbb{C}}$, the ordinary cup product $\cup$ is replaced by the quantum product $\cup_{q}$. 
Similar to the B-model, one can choose a symplectic basis $A , B$ of $H_{\mathrm{even}}(E, \mathbb{Z})$ to help describing the complexified K\"ahler structures.
Such a basis is determined by a marking
\begin{equation}\label{eqnAmodelmarking}
m: \mathbb{Z}\oplus \mathbb{Z}\cong H_{\mathrm{even}}(E,\mathbb{Z})\,.
\end{equation}
Then a generic local section $\Omega$ can be described in terms of the period integrals $\omega_{0}=\int_{A}\Omega\,, \omega_{1}=\int_{B}\Omega$ by
\begin{equation}
\Omega=\omega_{1}\beta+\omega_{0}\alpha\,,
\end{equation}
here $\alpha,\beta\in  H^{\mathrm{even}}(E,\mathbb{Z})$ is the dual basis of $A, B\in   H_{\mathrm{even}}(E,\mathbb{Z})$.

By comparing the datum in the A- and B-models, we can see that if we take the desired marking in \eqref{eqndesiredmarking} to be the one given in
\eqref{eqnAmodelmarking}, then indeed the corresponding $\mathrm{SL}_{2}(\mathbb{Z})$-action on $\mathcal{H}_{\mathrm{K}}$ 
meets all of the expectations from mirror symmetry.
The lattice $H_{\mathrm{even}}$ is usually referred to as the lattice of central charges in the language of stability conditions.\\

The above discussion also tells that the space $\mathcal{H}$ is actually not enough to capture the full information since the homothety is invisible on $\mathcal{H}$, which however is potentially useful in understanding the complete picture of mirror symmetry.
To illustrate this, we recall that in the B-model of $\check{E}$ to restore the homothety one should consider
the space of oriented basis in $\mathbb{R}^{2}$ given by $\{(\pi_{1},\pi_{0})| \mathrm{Im }(\pi_{1}/\pi_{0})>0\}$, see \cite{Katz:1976} for a review.
The $\mathrm{SL}_{2}(\mathbb{Z})$-action is described by
\begin{equation}
\gamma=
\begin{pmatrix}
a&b\\
c& d
\end{pmatrix}
\in
\mathrm{SL}_{2}(\mathbb{Z}): 
\begin{pmatrix}
\pi_{1}\\
\pi_{0}
\end{pmatrix}
\mapsto 
\begin{pmatrix}
a&b\\
c& d
\end{pmatrix}
\begin{pmatrix}
\pi_{1}\\
\pi_{0}
\end{pmatrix}
\,.
\end{equation}
If one only concentrates on the induced action on $\tau=\pi_{1}/\pi_{0}\in \mathcal{H}$, one 
only obtains the relative length of the vector $\Omega$ with respect to the marking $A,B$ and hence loses part of the information that $\Omega$ contains.

Similarly, in the A-model, the Chern character isomorphism gives
\begin{eqnarray}\label{eqnChern}
\mathrm{ch}: K_{0}(E)\otimes_{\mathbb{Z}}\mathbb{Q}&\rightarrow& H^{\mathrm{even}}(E,\mathbb{Z})\otimes_{\mathbb{Z}}\mathbb{Q}\,,\nonumber\\
\mathcal{E}&\mapsto& 
(c_{1}(\mathcal{E}),\mathrm{rank}(\mathcal{E}))\,.
\end{eqnarray}
If we specify a marking by taking $A,B$ to be the standard generators for $H_{0}(E,\mathbb{Z}),H_{2}(E,\mathbb{Z})$ respectively, then the
periods are $ (\omega_{1},\omega_{0})=(c_{1}(\mathcal{E}),\mathrm{rank}(\mathcal{E}))$.
Only keeping the slope of the pair $ (\omega_{1},\omega_{0})=(c_{1}(\mathcal{E}),\mathrm{rank}(\mathcal{E}))$
\begin{equation}\label{eqnslope}
\mu(\mathcal{E})={c_{1}(\mathcal{E})\over \mathrm{rank}(\mathcal{E})}\,,
\end{equation}
loses significant information of $\mathcal{E}$ and the $\mathrm{SL}_{2}(\mathbb{Z})$-action.

\subsection{Speculation: $S$-transform and Fourier-Mukai transform}

With the rank $2$ lattice $H_{\mathrm{even}}(E,\mathbb{Z})$ introduced into the story, the meaning of the $S$-transform is more clear. It acts by 
\begin{equation}
S=
\begin{pmatrix}
0&-1\\
1& 0
\end{pmatrix}
\in
\mathrm{SL}_{2}(\mathbb{Z}): (\omega_{1},\omega_{0})\mapsto  (\omega_{0},-\omega_{1})\,.
\end{equation}
Again by invoking the Chern character isomorphism, it amounts to saying that the sheaf $\mathcal{E}$ is sent to
$S(\mathcal{E})$ with 
\begin{equation}
c_{0}(S(\mathcal{E}))=-c_{1}(\mathcal{E})\,,\quad c_{1}(S(\mathcal{E}))=c_{0}(\mathcal{E})\,.
\end{equation}
The Chern character of the Abelian Fourier-Mukai transform of $\mathcal{E}\in \mathcal{D}^{b}(E)$, had $\mathcal{E}$ satisfied certain stability condition defined by the slope $\mu$ in 
\eqref{eqnslope}, satisfy exactly this equation. See \cite{Bartocci:2009fourier} for details.

This picture also seems to  be helpful in understanding conceptually the modularity in the enumeration of Hurwitz numbers of the elliptic curves.
It is natural to expect that there is a moduli space construction of the Hurwitz covers in terms of stable sheaves, and a generating series can be defined as 
a weighted sum. The Fourier-Mukai transform acts as an automorphism on this moduli space, and scales the weights by some overall factor.
Pulling out the overall factor gives the automorphy factor for the $S$-transform on the generating series as a modular form.
However, special care needs to be taken care of on the boundary of the moduli space which 
is expected to be responsible for the quasi-modularity rather than modularity. 
We are not able to make this fully rigorous so far and wish to pursue this line of thought in the future.

\section{Mirror manifolds of plane cubics and extra structures on lattices}
\label{sectorsion}

In the previous section we have already seen that the lattices play an important role in the mirror symmetry of the (universal) elliptic curve families.
This is one place where some arithmetic structures start to enter the story. 
We shall now discuss the mirror symmetry for plane cubic curves as another example.
In the sequel we shall see that structures like polarization and level structure determine the modular group symmetries on the moduli spaces on the two sides of mirror symmetry.

The construction for the mirror manifolds of CY varieties was firstly carried out in physics and known as the orbifold construction \cite{Greene:1990ud}. It was then realized as the polar duality between reflexive polytopes \cite{Batyrev:1994hm}. 
See \cite{Cox:2000vi} for a nice review of the history and also more detailed discussions.

\subsection{Toric duality}
\label{sectoricduality}

We shall first review the construction of the mirror manifolds of the plane cubics via toric duality following the toric language in the textbook \cite{Cox:2011}.

\subsubsection{Mirror of $\mathbb{P}^{2}$}

Consider the projective space $Y=\mathbb{P}^{2}$ as a toric variety. Its fan $\Sigma\subseteq N_{\mathbb{R}}$ is generated by the rays
\begin{equation}
\Sigma(1):\,\quad u_{1}=e_{1}\,,\quad u_{2}=e_{2}\,,\quad u_{3}=-e_{1}-e_{2}\,,\quad
\end{equation}
where
\begin{equation}
e_{1}=
\begin{pmatrix}
1\\
0
\end{pmatrix}\,,\quad
e_{2}=
\begin{pmatrix}
0\\
1
\end{pmatrix}\in N_{\mathbb{R}}\,.
\end{equation}
The corresponding reflexive polytope $\Delta$ is $3\Delta_{\mathrm{simp}}-(f_{1}+f_{2})=\mathrm{convex~hull}(v_{1},v_{2},v_{3})\subseteq M_{\mathrm{R}}$\,,
where
\begin{equation}
\quad v_{1}=2f_{1}-f_{2}\,,\quad v_{2}=-f_{1}+2f_{2}\,,\quad v_{3}=-f_{1}-f_{2}\,,
\end{equation}
and
\begin{equation}
f_{1}=
\begin{pmatrix}
1\\
0
\end{pmatrix}\,,\quad
f_{2}=
\begin{pmatrix}
0\\
1
\end{pmatrix}\in M_{\mathbb{R}}\,.
\end{equation}
See Figure \ref{figurepolyopefan} below.

\begin{figure}[h]
  \renewcommand{\arraystretch}{1.2} 
\begin{displaymath}
\includegraphics[scale=0.5]{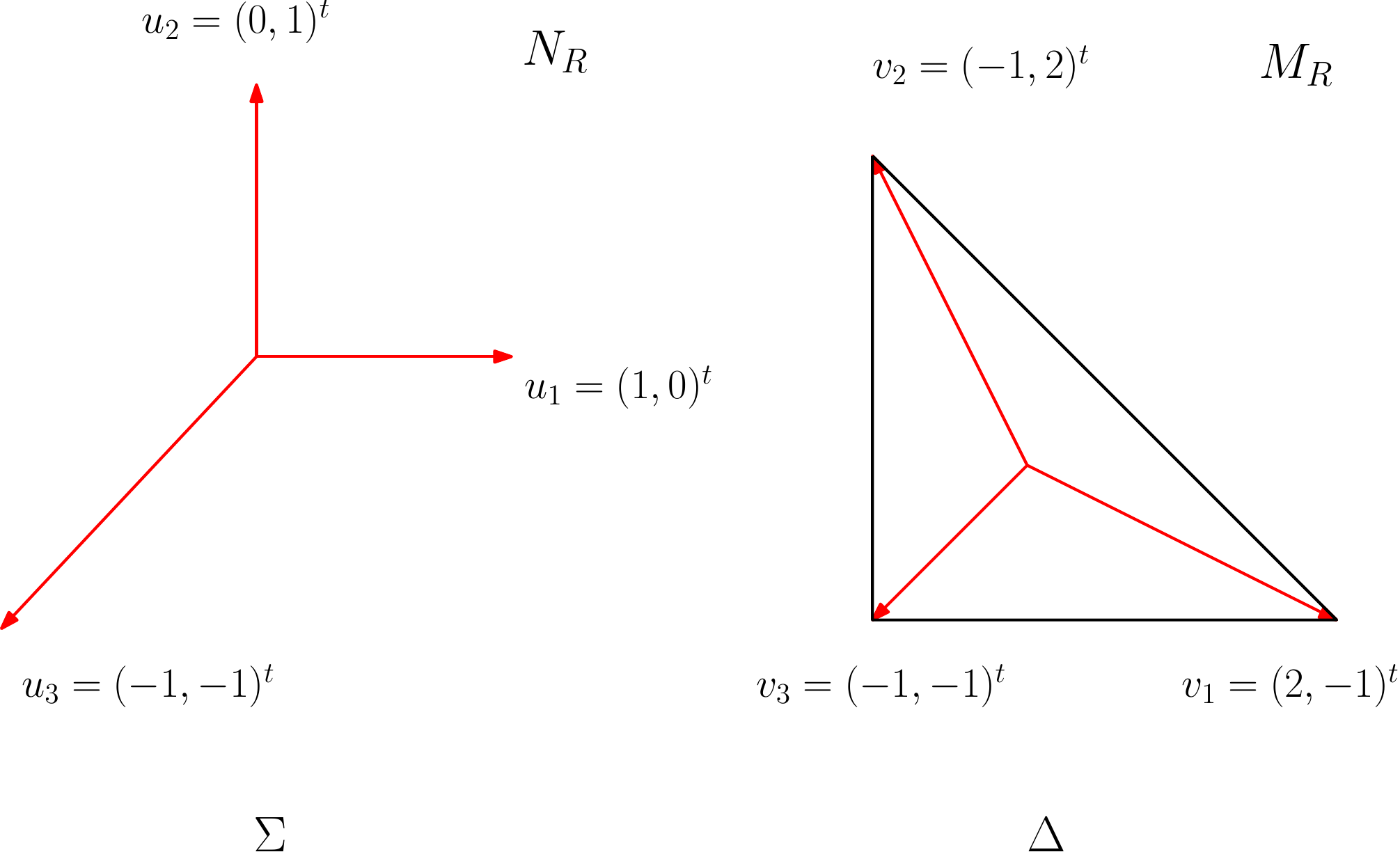}
\end{displaymath}
  \caption[polyopefan]{Fan and polytope of the toric variety $\mathbb{P}^{2}$.}
  \label{figurepolyopefan}
\end{figure}

The polar dual of the polytope $\check{\Delta} \subseteq \check{M}_{\mathbb{R}}=N_{\mathbb{R}} $ is the convex hull of the vertices $u_{1},u_{2},u_{3}$,
with the corresponding fan $\check{\Sigma} \subseteq \check{N}_{\mathbb{R}} =M_{\mathbb{R}}$ generated by $v_{1}, v_{2}, v_{3}$. 
This defines a new toric variety which is the mirror of $\mathbb{P}^{2}$, denoted by $\check{Y}=\check{\mathbb{P}^{2}}$ below. 
A more precise way to think about the new toric variety $\check{Y}$ is to regard $\check{\Sigma}$ as a stacky fan and then the corresponding variety as a toric stack. 

Now by standard construction, one has the homogeneous quotient description
\begin{equation}
Y=(\mathbb{C}^{3}-\{(0,0,0)\} )/ \mathbb{C}^{*}\,.
\end{equation}
One chooses the homogeneous toric coordinates on $\mathbb{C}^{3}$ to be $(z_{1},z_{2},z_{3})$, 
corresponding to the toric invariant divisors $D_{\rho},\,\rho=1,2,3$ respectively. Then the $\mathbb{C}^{*}$-action 
is given by
\begin{equation}
\lambda\in \mathbb{C}^{*}: (z_{1},z_{2},z_{3})\mapsto  (\lambda z_{1},\lambda  z_{2}, \lambda  z_{3})\,.
\end{equation}
Similarly, one has
\begin{equation}\label{eqnmirrorP2}
\check{Y}=(\mathbb{C}^{3}-\{(0,0,0)\} )/ (\mathbb{C}^{*}\times \mathbold{\mu}_{3})\,.
\end{equation}
Here $\mathbold{\mu}_{3}$ is the multiplicative cyclic group of order $3$.
By choosing the homogeneous toric coordinates on $\mathbb{C}^{3}$ to be $(x_{1},x_{2},x_{3})$, the action of $\mathbb{C}^{*}\times \mathbold{\mu}_{3}$ 
is then given by
\begin{equation}
(\lambda,\rho_{3})\in \mathbb{C}^{*}\times \mathbold{\mu}_{3}: (z_{1},z_{2},z_{3})\mapsto  (\lambda z_{1},\lambda  \rho_{3} z_{2}, \lambda  \rho_{3}^{2}z_{3})\,,\quad \rho_{3}=\exp({2\pi i\over 3})\,.
\end{equation}

\subsubsection{Mirror manifolds of plane cubics}

The toric duality \cite{Batyrev:1994hm} says that the mirror CY manifolds of the CY variety, represented as sections of $-K_{Y}$, are given by the sections of 
$-K_{\check{Y}}$.
A generic section of $-K_{Y}$ is determined by
\begin{equation}\label{eqnAmodelfamily}
\sum_{m\in \Delta\cap M} a_{m}\chi^{m}=0\,,\quad a_{m}\in\mathbb{C}\,,
\end{equation}
where $\{\chi^{m}\}$ are the character monomials corresponding to the lattice points $m\in \Delta\cap M$.
Switching to homogeneous toric coordinates $(z_{1},z_{2},z_{3})$ on $Y$, one has
\begin{equation}
\chi^{m}=\prod_{\rho\in \Sigma(1)} z_{\rho}^{\langle m,\, u_{\rho}\rangle+c_{\rho}}\,,
\end{equation}
where
the numbers $\{c_{\rho}\}$ are determined through the equation
\begin{equation}
-K_{Y}=\sum c_{\rho} D_{\rho}\,.
\end{equation}
In the present case, one has $c_{\rho}=1,\,\rho=1,2,3$.
Hence for example, we have
\begin{equation}
\chi^{v_{i}}=z_{i}^{3}\,,i=1,2,3\,, \quad \chi^{(0,0)}=z_{1}z_{2}z_{3}\,.
\end{equation}
The number of lattice points in $\Delta$ is $10$, corresponding to the number of cubic monomials in the homogeneous coordinates $(z_{1},z_{2},z_{3})$.
The equation in \eqref{eqnAmodelfamily} defines
a subvariety $\mathcal{X}$ in $\mathbb{C}^{3}_{(z_{1},z_{2},z_{3})}\times \mathbb{C}^{10}_{(a_{0},\cdots a_{9})}$.
Now we define a further $\mathbb{C}^{*}$-action given by
\begin{equation}
\mu\in \mathbb{C}^{*}: (a_{m})_{m\in \Delta\cap M}\mapsto  (\mu a_{m})_{m\in \Delta\cap M}\,.
\end{equation}
The quotient of the subvariety $\mathcal{X}$
by the product of the above two $\mathbb{C}^{*}$-actions
defines a family of cubic curves which we call the family of plane cubics.

The mirror family is given by 
\begin{equation}
 \sum_{\check{m}\in \check{\Delta}\cap \check{M}} b_{\check{m}}\chi^{\check{m}}=0\,, \quad b_{\check{m}}\in \mathbb{C}\,.
\end{equation}
There are only $4$ lattice points in $ \check{\Delta}\cap  \check{M}$, they give rise to the following monomials
\begin{equation}
\chi^{u_{i}}=x_{i}^{3}\,,i=1,2,3\,,  \quad \chi^{(0,0)}=x_{1}x_{2}x_{3}\,.
\end{equation}
Similar computations as above shows that the mirror family $\check{\mathcal{X}}$ is given by
\begin{equation}\label{eqnBmodelfamily}
\check{\mathcal{X}}: \{ b_{1}x_{1}^{3}+b_{2}x_{2}^{3}+b_{3}x_{3}^{3}+b_{0}x_{1}x_{2}x_{3}=0\}/ (\mathbb{C}^{*}\times \mathbold{\mu}_{3}\times \mathbb{C}^{*})\,.
\end{equation}
Here the $\mathbb{C}^{*}\times \mathbold{\mu}_{3}\times \mathbb{C}^{*}$-action is described by
\begin{equation}
(\lambda,\rho_{3},\mu)\in \mathbb{C}^{*}\times \mathbold{\mu}_{3}\times \mathbb{C}^{*}:
(x_{1},x_{2},x_{3},(b_{k})_{k=0}^{3})\mapsto 
(\lambda x_{1}, \lambda \rho_{3} x_{2}, \lambda \rho_{3}^{2}x_{3},(\mu b_{k})_{k=0}^{3})\,.
\end{equation}
Henceforward we shall ignore the $\mathbb{C}^{*}$-actions to simply the notations.

\subsubsection{Extra structure in the B-model}

For a generic elliptic curve, the dimensions of the space of complexified K\"ahler structures and of the complex structures are both one, it is hence a trivial statement that the dimensions involved on the two sides of mirror symmetry of plane cubics match.
However, the base of the elliptic curve family $\check{\mathcal{X}}$ given in \eqref{eqnBmodelfamily}
is not $\mathcal{H}_{\mathrm{C}}$ as this family is apparently different from the universal family of elliptic curves with the markings $\check{m}$.
It is therefore natural to ask what is the true moduli space and what is the extra structure carried by the members in the family.\\

To answer these questions, we recall that by scaling the homogeneous coordinates, which amounts to quotient by the $\mathrm{PGL}$-action and does not affect the discussions, the family $\check{\mathcal{X}}$ is equivalent to quotient of the so-called Hesse pencil $\mathcal{E}_{\mathrm{Hesse}}$ by the $\mathbold{\mu}_{3}$-action 
\begin{equation}\label{eqnquotientofHesse}
\{ x_{1}^{3}+x_{2}^{3}+x_{3}^{3}-3\psi x_{1}x_{2}x_{3}=0\}/ \mathbold{\mu}_{3}\,,
\end{equation}
where
\begin{equation}
\rho_{3}\in \mathbold{\mu}_{3}: (x_{1},x_{2},x_{3})\mapsto  (x_{1}, \rho_{3}x_{2}, \rho_{3}^{2}x_{3})\,.
\end{equation}
This gives the same result produced by the orbifold construction \cite{Greene:1990ud} for the mirror manifolds which will be reviewed below.

The Hesse pencil $\mathcal{E}_{\mathrm{Hesse}}$ is actually the universal family over the modular curve $\Gamma(3)\backslash\mathcal{H}^{*}$, where $\mathcal{H}^{*}$ is the compatification $\mathcal{H}\cup \mathbb{P}^{1}(\mathbb{Q})$.
This modular curve is the moduli space of pairs $(\check{E}, \check{m}_{3})$, where 
\begin{equation}\label{eqn3torsion}
\check{m}_{3}:\mathbb{Z}_{3}\oplus \mathbb{Z}_{3}\cong \check{E}[3]\,,
\end{equation}
with $\check{E}[3]$ being the group of $3$-torsion points of $\check{E}$.
This extra structure carried by the family is one kind of the level structure.
As a result, the family has $9$ flat sections corresponding to the $3$-torsion points (i.e., flex points) of the plane cubic curves. 
Computationally, that these sections are flat can be seen by observing that their projective coordinates are independent of the parameter $\psi$.
See \cite{Dolgachev:1997, Husemoller:2004, Artebani:2006, Dolgachev:2012} for detailed discussions.\\

Taking the invariants of the $\mathbold{\mu}_{3}$-action (without modulo projectivization) to be 
\begin{equation}\label{eqnetale}
X_{i}=x_{i}^{3}\,, \,i=1,2,3\,, \quad X_{0}=x_{1}x_{2}x_{3}\,,
\end{equation}
then the quotient is given by the following generically smooth complete intersection in $\mathbb{P}^{3}$
\begin{equation}
X_{1}+X_{2}+X_{3}-3\psi X_{0}=0\,, \quad X_{1}X_{2}X_{3}=X_{0}^{3}\,.
\end{equation} 
We remark that the above form for the mirror curve is what appears in the Hori-Vafa construction \cite{Hori:2000kt} for the mirror of some other closely related geometries.
By computing the canonical sheaf using the adjunction formula, we can see that indeed this is a CY variety of dimension one, hence an elliptic curve. Moreover,  by going to the affine coordinates 
\begin{equation}\label{eqnaffinecoordinatesofmirrorP2}
y=X_{2}/X_{3}\,,\quad x=-X_{0}/X_{3}\,,
\end{equation}
we are led to
\begin{equation}\label{eqnmirrorcubicequation}
y^{2}+3\psi x y+y=x^{3}\,.
\end{equation}
Now the family
\begin{equation}\label{eqnmirrorBfamily}
\check{\mathcal{X}}\cong \mathcal{E}_{\mathrm{Hesse}}/\langle \rho_{3}\rangle\,
\end{equation}
still has the modular curve $\Gamma(3)\backslash\mathcal{H}^{*}$ as its base and is also equipped with the $3$-torsion structure as discussed in \cite{Husemoller:2004}.
The two families $\mathcal{E}_{\mathrm{Hesse}}$ and $\check{\mathcal{X}}$ have equivalent variations of complex structures (excluding the $3$-torsion structures) and in particular the same Picard-Fuchs equations.
For this reason, in the literature when discussing some complex-geometric aspects of the B-model, one usually takes $\mathcal{E}_{\mathrm{Hesse}}$ as the mirror family. 

\begin{rem}
It is a classical result that the plane cubics with the $3$-torsion structure in \eqref{eqn3torsion} can be embedded 
equivariantly into $\mathbb{P}^{2}$ via theta functions.
Equivariance means that the translation actions by the elements in the group $\check{E}[3]$ in the domain
becomes some projective transformations on the image.
More precisely, with the particular embedding given in for example \cite{Dolgachev:1997}, one has, using $\check{E}\cong \mathbb{C}/(\mathbb{Z}\oplus \mathbb{Z}\tau)$ for some $\tau\in \mathcal{H}$\,,
\begin{eqnarray*}
{1\over 3}: (x_{1},x_{2},x_{3}) &\mapsto & (x_{1}, \rho_{3} x_{2}, \rho_{3}^{2}x_{3})\,,\\
{\tau\over 3}: (x_{1},x_{2},x_{3} )&\mapsto& (x_{2},x_{3}, x_{1})\,.
\end{eqnarray*}
Note that these translations do not preserve the origin of the elliptic curve.
The quotient by group of translations $\langle {1\over 3}\rangle$ is exactly the one that gives the above $3$-isogeny in \eqref{eqnmirrorBfamily}.
\qed
\end{rem}

\subsubsection{Speculation: extra structure in the A-model}

The philosophy of mirror symmetry implies that in the A-model of the plane cubics one should also see some extra structure.
The discussion in Section \ref{seclattice} suggests that the mirror should be the $3$-torsion in the lattice $H_{\mathrm{even}}(E,\mathbb{Z})$.

By going to the dual cohomology, this amounts to saying that there exist natural locally constant elements in $H^{\mathrm{even}}(E,\mathbb{Z})$ 
which generate a sub-lattice of index $9$, denoted symbolically by $3H^{\mathrm{even}}(E,\mathbb{Z})$.
Now under the Chern character isomorphism \eqref{eqnChern}, it suffices to find two sheaves $\mathcal{E}_{1}, \mathcal{E}_{2}$ on $E$ so that
$\mathrm{ch}(\mathcal{E}_{1}), \mathrm{ch}(\mathcal{E}_{2})$ generate $3H^{\mathrm{even}}(E,\mathbb{Z})$.
These sheaves should be locally constant along the moduli space direction. This would be the case if they are obtained by natural pull backs from the ambient $\mathbb{P}^{2}$ due to the simultaneous embedding of the fiber elliptic curves in the family. 
Now it is easy to see that taking $\mathcal{E}_{1}=(\mathcal{O}_{\mathbb{P}^{2}}^{\oplus 3})|_{E}, \mathcal{E}_{2}=(\mathcal{O}_{\mathbb{P}^{2}}(1))|_{E}$
does the job.

One can think of the elements in $H^{\mathrm{even}}$ as some generalized version of polarizations
which detect the "sizes" of the cycles in $H_{\mathrm{even}}$, then the above discussion implies that the torsion of polarization in the A-model is mirror to level structure in the B-model. This agrees with the constructions by SYZ mirror symmetry \cite{Strominger:1996mirror} or Fourier-Mukai transform \cite{Huybrechts:2006fourier, Bartocci:2009fourier}.

\subsection{Orbifold construction}

We mentioned in the above that the $3$-isogeny in \eqref{eqnmirrorBfamily} does not play a role in the complex geometry of the B-model.
At first glance, the statement sounds unlikely to be true since after all, it is $\check{\mathcal{X}}$ instead of $\mathcal{E}_{\mathrm{Hesse}}$ which is the mirror of the plane cubics.
The explanation is that although the quotient does not affect the complex geometry, it does has a non-trivial action on the K\"ahler and hence CY geometry of $\check{\mathcal{X}}$.
Looked from the mirror side, the mirror action does not affect the K\"ahler geometry of $\mathcal{X}$ much but is important for the complex geometry of $\mathcal{X}$.

\subsubsection{Orbifold construction for the mirror quintic family}

Again to motivate the discussion, we first recall the orbifold construction \cite{Greene:1990ud} for the mirror of the quintic family, following the exposition in \cite{Gross:2003calabi}. 

Consider a generic quintic $Q$ in $\mathbb{P}^{4}$.
Straightforward computation shows that 
\begin{equation}
h^{1,1}(Q)=1\,,\quad h^{2,1}(Q)=101\,.
\end{equation}
The space of K\"ahler structures is spanned by the pull back of the hyperplane class in the ambient space $\mathbb{P}^{4}$. 
For the complex structure deformations, heuristically they correspond to the vector space $V$ of degree-$5$ monomials in the $5$ homogeneous coordinates on $\mathbb{P}^{4}$.
The group $\mathrm{PGL}_{5}(\mathbb{C})$ acts by bi-holomorphisms. After projectivization, the space of complex structure deformations has dimension
\begin{equation}
{5+5-1\choose 5-1}-\mathrm{dim}\mathrm{PGL}_{5}(\mathbb{C})-1=101\,.
\end{equation}
This then yields the universal family $\mathcal{Q}$ of quintics in $\mathbb{P}^{4}$.

Now mirror symmetry predicts that a generic member $\check{Q}$ in the mirror family  $\check{\mathcal{Q}}$ should satisfy
\begin{equation}
h^{1,1}(\check{Q})=101\,,\quad h^{2,1}(\check{Q})=1\,.
\end{equation}
Consider the action of diagonal symmetries 
\begin{equation}
G=\mathbb{P}\{(\rho_{5}^{n_{1}},\rho_{5}^{n_{2}},\rho_{5}^{n_{3}},\rho_{5}^{n_{4}},\rho_{5}^{n_{5}})\in 
\mathbold{\mu}_{5}^{5}~|\sum_{i=1}^{5} n_{i}\equiv 0~\mathrm{mod}\,5\}\subseteq 
\mathrm{PGL}_{5}(\mathbb{C})\,, \,\,\rho_{5}=\exp({2\pi i\over 5})\,.
\end{equation}
One can understand this group 
as the one in the homogeneous quotient in constructing the mirror $\check{\mathbb{P}}^{4}=\mathbb{P}^{4}/G$
of $\mathbb{P}^{4}$ in a way similar to 
\eqref{eqnmirrorP2} using toric duality.
The vector space $V$ then decomposes into sums of representations
\begin{equation}
V=V_{0}\oplus \bigoplus_{\chi} V_{\chi}\,.
\end{equation} 
Here in the direct sum $\chi$ runs over the space of non-trivial characters.
The representation $V_{0}$ with trivial character is spanned by
\begin{equation}
V_{0}=\mathbb{C}[x_{1}^{5},  x_{2}^{5},x_{3}^{5},x_{4}^{5},x_{5}^{5}, x_{1}x_{2}x_{3}x_{4}x_{5}]\,.
\end{equation}
Under the action of $\mathrm{PGL}_{5}(\mathbb{C})$, the space $V_{0}$ gives an one-parameter family of quintics in $\mathbb{P}^{4}$, called the Dwork pencil,
\begin{equation}
\mathcal{Q}_{\mathrm{Dwork}}:\quad \sum_{i=1}^{5} x_{i}^{5}-5\psi \prod_{i=1}^{5}x_{i}=0\,.
\end{equation}
Then one takes the quotient $\mathcal{Q}_{\mathrm{Dwork}}/G$ of $\mathcal{Q}_{\mathrm{Dwork}}$ by the $G$-action. The members in the resulting family are singular varieties.
One applies a resolution of singularity to obtain a family of smooth varieties.
The mirror manifold is declared to be the resolution
\begin{equation}
\check{\mathcal{Q}}=\widetilde{\mathcal{Q}_{\mathrm{Dwork}}/G}\,.
\end{equation}
Note that the members $\check{Q}$ in this family are not simultaneously embedded into $\mathbb{P}^{4}$ anymore due to the 
quotient which produces singularities and brings in the resolution.
Again a heuristic counting of $h^{1,1}(\check{Q})$ is given as follows.
The resolution creates exceptional divisors which under the Poincar\'e dual contribute to $h^{1,1}(\check{Q})$, then
\begin{equation}
h^{1,1}(\check{Q})=h^{1,1}(Q)+\mathrm{number ~of ~exceptional~divisors}\,.
\end{equation}
The counting gives the desired number $101$. 
See \cite{Gross:2003calabi} for references which offer rigorous approaches in proving the results on dimensions.

\subsubsection{Speculation: roles of group action}

In retrospect, we can see that the group $G$ plays two roles.
The first is to cut down the dimension of space of complex structure deformations by keeping only the $G$-invariants.
The resulting family is still simultaneously embedded in $\mathbb{P}^{4}$.
The second is to increase the dimension of space of K\"ahler structure deformations by forming the quotient followed by resolution, at the price of losing the simultaneous embedding.

If we only care about the variation of complex structures of $\check{\mathcal{Q}}$, for example if we only focus on variation of periods, then
$\mathcal{Q}_{\mathrm{Dwork}}$ is already enough.
This is in fact what is customarily done in the literature: 
one uses
the Picard-Fuchs equation for $\mathcal{Q}_{\mathrm{Dwork}}$
in computing the periods of $\check{\mathcal{Q}}$, see \cite{Candelas:1990rm, Gross:2003calabi}.
Note also that if we regard $\mathcal{Q}_{\mathrm{Dwork}}/G$ as an orbifold, then the additional K\"ahler structure deformations 
coming from the resolution correspond to
some twisted sectors in the Chen-Ruan cohomology \cite{Chen:2001} of the orbifold. Hence the orbifold and its resolution contain the same amount of information on the K\"ahler structure deformations.

\begin{figure}[h]
  \renewcommand{\arraystretch}{1.2} 
\begin{displaymath}
\xymatrix{
\mathcal{Q} \ar[ddrr]^{G} \ar[rrrr]^{\mathrm{mirror}} &&&& \check{\mathcal{Q}}=\mathcal{Q}_{\mathrm{Dwork} }/G\\
&&&&\\
       && \mathcal{Q}_{\mathrm{Dwork}}  \ar[uurr]^{G}   & & 
}
\end{displaymath}
  \caption[quintic]{Transition in the deformation space of CY structures for the quintic}
  \label{figurequintic}
\end{figure}
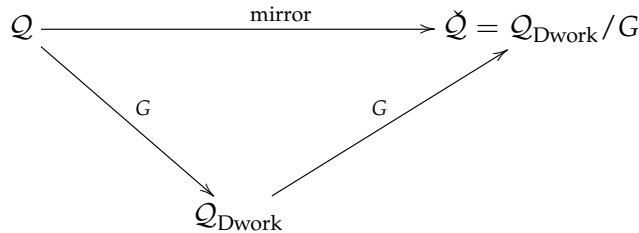

Now we can interpret the above results in terms of the diagram in Figure \ref{figurequintic}.
The $G$-action on the down-right arrow is used to pick out the $G$-invariants. 
It does not change the dimension of the space of K\"ahler structure deformations.
The inverse changes the dimension of the space of complex structure deformations by turning on the deformations along the monomials in the space $\bigoplus_{\chi}V_{\chi}$.
While the $G$-action on the up-right arrow is used to form the orbifold quotient, it does not change the dimension of complex structure deformations.
Its inverse keeps the $G$-invariants of the $H^{2}$ part in the Chen-Ruan cohomology of the orbifold $\mathcal{Q}_{\mathrm{Dwork} }/G$.

\subsubsection{Speculation: back to plane cubics}

Now we extend what we have learned from the quintic family case to the plane cubics. 
Note the dimensions of the spaces of K\"ahler structures and complex structures for the elliptic curve and its mirror, or the corresponding Hodge numbers, are always one.\
In order to figure out the different roles of the families that are involved in the mirror construction, 
what we should be focusing on are the dimensions of the space of deformations created by the torsions in the corresponding lattices, which we call twisted sectors borrowing the terminology from physics.
We use the notations $(\mathrm{tdim}_{\mathrm{Kahler}}, \mathrm{tdim}_{\mathrm{complex}})$ to denote the corresponding dimensions.

\begin{figure}[h]
  \renewcommand{\arraystretch}{1.2} 
\begin{displaymath}
\xymatrix{
\mathcal{X}:  (9,1) \ar[ddrr]^{G} \ar[rrrr]^{\mathrm{mirror}} &&&&  \check{\mathcal{X}}=\check{\mathcal{E}}_{\mathrm{Hesse} }/G: (1,\check{9})\\
&&&&\\
       && \mathcal{E}_{\mathrm{Hesse}}: (9,9)  \ar[uurr]^{G}   & & 
}
\end{displaymath}
  \caption[cubic]{Transition in the deformation space of CY structures for the cubic}
  \label{figurecubic}
\end{figure}
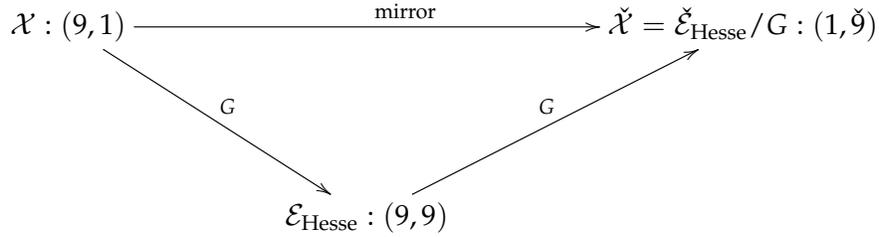

By analogy from the quintic family we are then led to the diagram in Figure \ref{figurecubic}.

We now interpret the earlier results discussed in Section \ref{sectoricduality} using this diagram. 
For the family $\mathcal{X}$, one has $\mathrm{tdim}_{\mathrm{Kahler}}=9$ since the lattice $H_{\mathrm{even}}(E,\mathbb{Z})$ is decorated with its $3$-torsion. The statement that $\mathrm{tdim}_{\mathrm{complex}}=1$ comes from the dimension count of the vector space of cubic monomials modulo $\mathrm{PGL}$-action.

The $G$-action on the down-right arrow keeps only the $G$-invariant cubic monomials, leaving
the $3$-torsion in $H_{\mathrm{even}}(E,\mathbb{Z})$ unchanged as the simultaneous embedding into $\mathbb{P}^{2}$ remains. In addition, the resulting family $\mathcal{E}_{\mathrm{Hesse}}$ has the $3$-torsion structure in $H_{1}(E,\mathbb{Z})$.
For the up-right arrow, the orbifold quotient on the family $\mathcal{\mathcal{E}}_{\mathrm{Hesse}}$
does not preserve the original simultaneous embedding into $\mathbb{P}^{2}$ and hence the $3$-torsion in $H_{\mathrm{even}}(E,\mathbb{Z})$ is lost.
This can be seen from the fact that the map in \eqref{eqnetale} destroys the original ambient space $\mathbb{P}^{2}$.
This results in the counting $\mathrm{tdim}_{\mathrm{Kahler}}=1$.
However, as explained in the paragraph below \eqref{eqnmirrorBfamily}, the new family carries a new
$3$-torsion structure in $H_{1}$ and thus $\mathrm{tdim}_{\mathrm{complex}}=9$.
The notation $\check{9}$ means that this new $3$-torsion structure is different from the one on $\mathcal{\mathcal{E}}_{\mathrm{Hesse}}$\,.\\

It appears that studying the torsion structures offer some new insights in understanding the modularity in the Gromov-Witten theory of the elliptic orbifold curves, see \cite{Shen:2014, Lau:2014, Shen:2016}. These structures, which are rooted in the symmetries of the families, are also closely related to the 
studies of oscillating integrals and GKZ hypergeometric functions in the appearance of chain integrals. This will be discussed in a forthcoming work.

\section{Detecting dualities from Picard-Fuchs equations}
\label{secdualities}

For either the family $\mathcal{E}_{\mathrm{Hesse}}$ or $\check{\mathcal{X}}$, 
the fibers over $\psi$, $\rho_{3}\psi$ and $\rho_{3}^{2}\psi$ are isomorphic. 
In particular, for the Hesse pencil, the isomorphism is induced by the projective transformation
\begin{equation}
[x_{1},x_{2},x_{3}]\mapsto [\rho_{3} x_{1},x_{2},x_{3}]\,.
\end{equation}
For simplification, one usually applies a base change $\psi\mapsto \alpha=\psi^{-3}$.
The resulting two families would then have the base being a copy of $\mathbb{P}^{1}$ parametrized by $\alpha$. They are again related by the $3$-isogeny mentioned before and have equivalent variations of complex structures.

For definiteness, we shall discuss the base change of $\check{\mathcal{X}}$, which we still denote by the same symbol by abuse of notation, 
\begin{equation}\label{eqnmirrorcubicbasechange}
\check{\mathcal{X}}: \quad y^{2}+3x y+y=\alpha x^{3}\,.
\end{equation}
It is known that $\psi,\alpha$ are the Hauptmoduln for the modular groups $\Gamma(3),\Gamma_{0}(3)$, respectively.  Moreover, the family $\check{\mathcal{X}}$ in
\eqref{eqnmirrorcubicbasechange} is actually the universal elliptic curve family over the modular curve $\Gamma_{0}(3)\backslash \mathcal{H}^{*}$
which is the moduli space of pairs
$(C, H)$, where $C$ is an elliptic curve and $H$ is a cyclic subgroup of order $3$
of $C[3]$.
The underlying modularity and moduli space interpretation is very useful in analyzing symmetries of the families, see for example \cite{Candelas:1993, Artebani:2006}.

\subsection{Fricke involution}
There is a particular involution on the moduli space $\Gamma_{0}(3)\backslash \mathcal{H}^{*}$ given by
\begin{equation}
(C,H)\mapsto (C/H, C[3]/H)\,.
\end{equation}
It is represented by the Fricke involution on $\mathcal{H}^{*}$
\begin{equation}
W: \tau\mapsto -{1\over 3\tau}\,.
\end{equation}
With a suitable choice for the Haupmodul $\alpha(\tau)$, see \cite{Maier:2009}, the Fricke involution
induces the affine transformation 
\begin{equation}
\alpha(\tau)\mapsto \beta(\tau):=1-\alpha(\tau)\,.
\end{equation}
Note that although the two points $\alpha=0,1$, corresponding to $\tau=i\infty, 0$ have the same $j$-invariant, the above transformation that is induced from the moduli interpretation is not the $S$-transform
but the Fricke involution $W$ which does not lie in $\mathrm{SL}_{2}(\mathbb{Z})$. In particular it does not belong to the Deck group of the covering $ \alpha\mapsto j(\alpha)$.
This involution plays a very important role in the mirror symmetry of some special non-compact CY threefolds \cite{Alim:2013eja}.
Also in the Seiberg-Witten theory \cite{Seiberg:1994rs}, the Fricke involution is what underlies the electro-magnetic duality, see 
\cite{Zhou:2014thesis} for more discussion on this.\\

One could have found this involution on the family without using the underlying modularity
but by using the Picard-Fuchs equation for this family. 
In the studies of mirror symmetry for more general CY varieties,
the former is usually difficult to make sense of while
the latter is typically what one can easily obtain. 
We now explain how this works following the discussions in \cite{Alim:2013eja, Zhou:2014thesis}.\\

Recall that the Picard-Fuchs operator for the family $\check{\mathcal{X}}$ in \eqref{eqnmirrorcubicbasechange} is 
\begin{equation}
\mathcal{L}_{\mathrm{PF}}=\theta_{\alpha}^{2}-\alpha(\theta_{\alpha}+{1\over 3})(\theta_{\alpha}+{2\over 3})\,,\quad \theta_{\alpha}:=\alpha{\partial\over \partial \alpha}\,.
\end{equation}
This can be deduced by using, for example, the Griffiths-Dwork method for the Hesse pencil
and then use the fact that the isogeny does not affect the Picard-Fuchs equation.

The resulting Picard-Fuchs equation $\mathcal{L}_{\mathrm{PF}}\pi=0$ has three singularities located at 
$\alpha=0,1,\infty$. 
Around each point, one chooses an appropriate local coordinate to rewrite the Picard-Fuchs equation and then can find a local basis in terms of hypergeometric functions.
In particular, near the point $\alpha=1$, one adapts the local coordinate
$\beta=1-\alpha$ and rewrite the Picard-Fuchs equation as
\begin{equation}
\mathcal{L}_{\mathrm{PF}}\pi=\left(\theta_{\beta}^{2}-\beta(\theta_{\beta}+{1\over 3})(\theta_{\beta}+{2\over 3})\right)\pi=0\,,\quad \theta_{\beta}:=\beta{\partial\over \partial \beta}\,.
\end{equation}
The Picard-Fuchs operator in the $\beta$-coordinate takes exactly the same form as the one written
in the $\alpha$-coordinate. This means that if $\pi(\alpha)$ is a solution, then $\pi(\beta)$
is a solution as well. Hence one finds the symmetry $\alpha\mapsto \beta=1-\alpha$ on the level of periods.

Moreover, one can choose a suitable basis $\pi_{0}(\alpha),\pi_{1}(\alpha)$ of solution near $\alpha=0$
and then define the normalized period to be
\begin{equation}\label{eqntau}
\tau(\alpha)={\pi_{1}(\alpha)\over \pi_{0}(\alpha)}\,.
\end{equation}
It has logarithm growth as $\alpha$ goes to zero and lies on the upper-half plane.
Moreover, a particular basis $\pi_{0}(\alpha)=\,_{2}F_{1}({1\over 3},{2\over 3};1;\alpha),\pi_{1}(\alpha)={i\over \sqrt{3}}\,_{2}F_{1}({1\over 3},{2\over 3};1;\beta)$ can be chosen such that the local monodromy near $\alpha=0$ is given by $\tau(\alpha)\mapsto \tau(\alpha)+1$. These two solutions are period integrals over the
vanishing cycles $A,B$ at the singularities $\alpha=0,1$ respectively.
\begin{rem}
In fact, this choice of basis satisfies the property \cite{Berndt:1995} that $\tau(\alpha)=\tau$, where the parameter $\tau$ is the modular variable satisfying 
\begin{equation}
j(\tau)={1\over q}+744+\cdots\,,\quad q=e^{2\pi i \tau}\,.
\end{equation}
\qed
\end{rem}
Then one can check that
\begin{equation}
\tau(\beta)=-{1\over 3\tau(\alpha)}\,.
\end{equation}
Hence one recovers the Fricke involution.
That is, the monodromy consideration naturally singles out the correct variable $\tau(\alpha)$ and the correct form for the involution on the moduli space.
 
\subsection{Cayley transform}

Similarly, one can take a suitable basis of periods near the orbifold point $\alpha=\infty$
and then form the normalized period $\tau_{\mathrm{orb}}$ there. By using the analytic continuation formulas of hypergeometric series \cite{Erdelyi:1981}, it follows that
\begin{equation}
\tau_{\mathrm{orb}}={\tau(\alpha)-\tau_{*}\over \tau(\alpha)-\overline{\tau_{*}}}\,,
\end{equation}
for some $\tau_{*}$ with $\mathrm{Im}\tau_{*}>0$.
This factional linear transform between the normalized periods is the Cayley transform based at the point $\tau_{*}$ on the upper-half plane.
It turns out that this simple transformation, induced by analyzing the Picard-Fuchs equation, is \cite{Shen:2016} what underlies the Calabi-Yau/Landau-Ginzburg correspondence \cite{Witten:1993} for the elliptic orbifold curves.
This correspondence is another interesting physics duality which has attracted a lot of attention in mathematics.\\

The idea of using the Picard-Fuchs equation to detect dualities 
 applies to other models of elliptic curve families and also some related geometries, see e.g., \cite{Alim:2013eja, Shen:2016} and references therein.


\section{Picard-Fuchs equations and Yukawa couplings}
\label{secYukawa}

One of the important predictions of mirror symmetry is that calculation on the genus zero Gromov-Witten invariants of
a CY variety can be turned into computation on period integrals of the mirror CY variety \cite{Candelas:1990rm}.
In this section, we shall review how this works by discussing a few examples. In the course we shall also see that the Picard-Fuchs equations
can be very useful in studying the Weil-Petersson metric on the deformation space of CY varieties.

\subsection{Elliptic curves}

It is well-known that the genus zero Gromov-Witten invariants of the elliptic curve $E$ are trivial in nonzero degree. One can see this by employing the definition of these invariants and checking that the integrals are always trivial due to dimension reasons \cite{Cox:2000vi}. 
This then implies the following identity on the generating series
of these invariants called Yukawa coupling,
\begin{equation}
C_{t}:=\sum_{d\geq 0}N_{0,d}\,q_{t}^{d}=1\,,\quad q=\exp(2\pi \sqrt{-1} t)\,.
\end{equation}
The numbers $\{N_{0,d}\}$ are the genus zero degree $d$ Gromov-Witten invariants\footnote{More precisely, what is discussed here is the generating series of Gromov-Witten invariants of genus zero, with one marking. Similarly, for the K3 surfaces and CY threefolds discussed below, the Yukawa couplings are those with two and three markings, respectively.}.

The miracle of mirror symmetry says that \cite{Candelas:1990rm} the above generating series can be computed in the B-model through
\begin{equation}
C_{\tau}:={1\over (\int_{A}\Omega)^{2}}\int_{\check{E}}\Omega\wedge \partial_{\tau}\Omega\,.
\end{equation}
Here the cycle $A$ is the vanishing cycle near $\tau=\sqrt{-1}\infty$ and
the holomorphic top form $\Omega$ is given in \eqref{eqnBtopform}.
The mirror map sends the K\"ahler structure parameter $t$ to the complex structure $\tau$, as discussed in Section \ref{secintro}.
An easy computation shows that indeed that 
\begin{equation}\label{eqnYukawauniversal}
C_{\tau}=1\,.
\end{equation}

Since the Gromov-Witten invariants are deformation invariant, in the A-model we can take any reasonably behaved family. 
Then by mirror symmetry, the family in the B-model could be any nicely behaved family.
In particular we can take the plane cubics to be the A-model.
Then the mirror is the family $\check{\mathcal{X}}$ described earlier in \eqref{eqnmirrorcubicbasechange}.
Now the Yukawa coupling is given by
\begin{equation}\label{eqnYukawamodular}
C_{\tau}={1\over \pi_{0}^{2}}{1\over 2\pi i}{\partial \alpha\over \partial \tau}C_{\alpha}\,,
\end{equation}
with
 \begin{equation}\label{eqnYukawaalgebraic}
C_{\alpha}=\int_{\check{\mathcal{X}}_{\alpha}} \Omega\wedge \partial_{\alpha}\Omega\,.
\end{equation}
Here $\pi_{0}$ is the integral of the holomorphic top form $\Omega$ on the vanishing cycle $A$ near the point $\alpha=0$, coresponding to $\tau=\sqrt{-1}\infty$
according to \eqref{eqntau}.
The quantity $C_{\alpha}$ satisfies a first order differential equation which is easily derived from the Picard-Fuchs equation
$\mathcal{L}_{\mathrm{PF}}\pi=0$. 
Solving this equation, one gets a rational function
\begin{equation}
C_{\alpha}={1\over \alpha(1-\alpha)}\,.
\end{equation}
Then the result $C_{\tau}=1$ follows from the Schwarzian equation for $\tau(\alpha)$. See \cite{Zhou:2013hpa} for details.

Now if one applies the same discussion to a general algebraic family of elliptic curves, then by using \eqref{eqnYukawauniversal}, \eqref{eqnYukawamodular} and \eqref{eqnYukawaalgebraic}, one
would produce interesting identities which are otherwise difficult to check directly.
For example, consider the $E_{8}$ elliptic curve family, see \cite{Lian:1994zv},
\begin{equation}
x_{1}^{6}+x_{2}^{3}+x_{3}^{2}-({\alpha\over 432})^{-{1\over 6}}x_{1}x_{2}x_{3}=0\, , \quad
j(\alpha)={1\over \alpha(1-\alpha)}\,.
\end{equation}
One gets the following identities, see also \cite{Hosono:2008ve, Zhou:2013hpa},
\begin{equation}
{1\over 2\pi i}\partial_{\tau}\alpha(\tau)=j(\tau)^{-1}\,_{2}F_{1}({1\over 6},{5\over 6};1;\alpha(\tau))^{2}\,,
\quad j(\tau)={1\over \alpha(\tau)(1-\alpha(\tau))}\,.
\end{equation}
\begin{rem}
Alternatively, all of these can be checked by using the fact that the parameter $\alpha(\tau)$ is the Haupmodul for some modular group and has
very nice expressions in terms of $\eta$-or $\theta$-functions, and that the periods are also related to modular forms. See \cite{Zhou:2013hpa}
for detailed discussions.
\qed
\end{rem}

\subsection{K3 surfaces}

Similar to the elliptic curves, the Gromov-Witten invariants of K3 surfaces are also trivial. Then the same reasoning above is supposed to yield 
non-trivial identities involving special functions. See \cite{Nagura:1993kd, Lian:1994zv, Lian:1995, Lian:1996} for related discussions.

We consider the Dwork pencil as the B-model for example. The equation for the family is given by
\begin{equation}
x_{1}^{4}+x_{2}^{4}+x_{3}^{4}+x_{4}^{4}-4z^{-{1\over 4}}x_{1}x_{2}x_{3}x_{4}=0\,,\quad z\in \mathbb{P}^{1}\,.
\end{equation}
The Picard-Fuchs operator is given by the hypergeometric differential operator
\begin{equation}
\mathcal{L}_{\mathrm{K3}}=\theta_{z}^{3}-z(\theta_{z}+{1\over 4})(\theta_{z}+{2\over 4})(\theta_{z}+{3\over 4})\,,\quad  \theta_{z}=z{\partial\over \partial z}\,.
\end{equation}
An easy computation \cite{Cox:2000vi} gives the Yukawa coupling in the parameter $z$
\begin{equation}
C_{zz}={1\over z^{2}(1-z)}\,.
\end{equation}
By considering the indicial equation at the regular singular point $z=0$, we know that we can choose a basis so that only one of them is regular, the other two have $\log z,(\log z)^{2}$ behavior near $z=0$.
We denote these three periods near $z=0$ by $\pi_{i},i=0,1,2$, respectively.
For example, we can take them to be the ones obtained by using the Frobenius method
\begin{eqnarray*}
\pi_{0}(z)&=&\pi_{0}(z,\rho)|_{\rho=0}=\,_{3}F_{2}({1\over 4},{2\over 4},{3\over 4};1,1; z)\,,\\
\pi_{1}(z)&=&c_{1}\partial_{\rho}\pi_{0}(z,\rho)|_{\rho=0}\\
&=&c_{1}\left(\pi_{0}(z)\ln ({z\over 4^{4}})+\sum_{n\geq 1} {\Gamma(4n+1)\over \Gamma(n+1)^{4} } (4\psi(4n+1)-4\psi(n+1)) ({z\over 4^{4}})^{n}\right)\,,\\
\pi_{2}(z)&=&c_{2}\partial_{\rho}^2\pi_{0}(z,\rho)|_{\rho=0}\,.
\end{eqnarray*}
Here 
\begin{equation}
\pi_{0}(z,\rho)=\sum_{n\geq 0} {\Gamma(1) \Gamma(1)\over \Gamma({1\over 4}) \Gamma({2\over 4}) \Gamma({3\over 4})}{\Gamma({1\over 4}+n+\rho)  \Gamma({2\over 4}+n+\rho) \Gamma({3\over 4}+n+\rho)
\over \Gamma(1+n+\rho) \Gamma(1+n+\rho)  } {z^{n+\rho}\over \Gamma(1+n+\rho)}\,,
\end{equation}
and $\psi(z)=\partial_{z}\ln\Gamma(z)$, while $c_{1}, c_{2}$ are some constants.
We also define the normalized period to be 
\begin{equation}
s={\pi_{1}\over \pi_{0}}\,.
\end{equation}
Note that a different choice for the basis $(\pi_{0},\pi_{1},\pi_{2})$, with the constraints that for some constants $c_{1},c_{2}$
\begin{equation}\label{eqnconstraints}
\pi_{0}=\mathrm{regular}\,,\quad \pi_{1}\sim c_{1}\log z+\cdots\,,\quad  \pi_{2}\sim c_{2}(\log z)^2+\cdots\,,
\end{equation}
induces an affine transformation on $s$ and hence a scaling on the Yukawa coupling
\begin{equation}
C_{ss}={1\over \pi_{0}^{2}} {1\over (2\pi i)^2}({\partial z\over \partial s})^{2}C_{zz}\,.
\end{equation}
It is known that up to scaling the parameter $s$ is mirror to the K\"ahler structure parameter $t$ in the A-model, see \cite{Cox:2000vi}.
Since $C_{tt}=1$, we know
\begin{equation}
C_{ss}=c\,,
\end{equation}
for some constant $c$. Therefore we are led to
\begin{equation}\label{eqnSchwarzianK3}
{1\over (2\pi i)^{2}}({\partial z\over \partial s})^{2}=c \pi_{0}^{2} z^{2}(1-z)\,.
\end{equation}
This then gives a non-trivial identity involving the hypergeometric series $\pi_{1},\pi_{0}$.\\

Having carried out the computations using Picard-Fuchs equation, one might wonder whether there exists some structure like modularity which underlies this, similar to the elliptic curve case. The answer is affirmative.
To proceed, we first note that the Picard-Fuchs operator $\mathcal{L}_{\mathrm{K3}}$ is actually is the symmetric square of some second order differential operator
\cite{Nagura:1993kd, Lian:1994zv} 
\begin{equation}\label{eqntriangular}
\mathcal{L}_{\mathrm{triangular}}=\theta_{z}^{2}-z(\theta_{z}+{1\over 8})(\theta_{z}+{3\over 8})\,,
\end{equation}
in the sense that 
\begin{equation}\label{eqsymmetricsquare}
\mathrm{Solution}(\mathcal{L}_{\mathrm{K3}})=\mathrm{Sym}^{\otimes 2}( \mathrm{Solution}(\mathcal{L}_{\mathrm{triangular}}))\,.
\end{equation}
Near $z=0$ one can take a basis of solutions to $\mathcal{L}_{\mathrm{triangular}}$ to be 
\begin{equation}
u_{0}(z)=\,_{2}F_{1}({1\over 8}, {3\over 8};1;z)\,, \quad u_{1}(z)=\,_{2}F_{1}({1\over 8}, {3\over 8};{1\over 2};1-z)\,.
\end{equation}
The second solution $u_{2}$ is interpreted as the analytic continuation near $z=0$ and hence has $\log$ behavior.
One can then check the symmetric square structure directly. 
For example, the relation $\pi_{0}=u_{0}^{2}$ follows from Clausen's identity.

Now we pick the following basis satisfying the conditions in \eqref{eqnconstraints}
\begin{equation}
\pi_{0}=u_{0}^{2}\,,\quad  \pi_{1}=u_{0}u_{1}\,,\quad  \pi_{2}=u_{1}^{2}\,.
\end{equation}
Then we get
\begin{equation}
s={\pi_{1}\over \pi_{0}}={u_{1}\over u_{0}}\,.
\end{equation}
We then treat $s$ as the normalized period for the equation $\mathcal{L}_{\mathrm{triangular}}u=0$.
This gives the relation \cite{Erdelyi:1981}
\begin{equation}
{\partial z\over \partial s}=z(1-z)^{1\over 2}u_{0}^{2}\,.
\end{equation}
Then \eqref{eqnSchwarzianK3} follows easily.

\begin{rem}
The above discussion only used the symmetric square structure and no relation to elliptic curve family or modular forms is relied on.
In fact, the parameter $z$ is the Hauptmodul for the modular group $\Gamma_{0}^{+}(2)$, where the triangular group $\Gamma_{0}^{+}(2)$
is the Fricke extension of the modular group $\Gamma_{0}(2)<\mathrm{SL}_{2}(\mathbb{Z})$.
The corresponding "universal" (modulo the issue of orbifold) elliptic curve family over the modular curve $\Gamma_{0}(2)\backslash\mathcal{H}^{*}$ is the $E_{7}$ elliptic curve family
\begin{equation}\label{eqnE7family}
x_{1}^{4}+x_{2}^{4}+x_{3}^{2}-({\alpha\over 64})^{-{1\over4}}x_{1}x_{2}x_{3}=0\,.
\end{equation}
The Picard-Fuchs operator is 
\begin{equation}
\mathcal{L}_{\mathrm{elliptic}}=\theta_{\alpha}^{2}-\alpha (\theta_{\alpha}+{1\over 4}) (\theta_{\alpha}+{3\over 4})\,,
\end{equation}
where $\alpha$ is the Haupmodul for $\Gamma_{0}(2)$ given in e.g., \cite{Maier:2009}.
The Haupmodul $z$ is related to $\alpha$ by
\begin{equation}
z=4\alpha(1-\alpha)\,.
\end{equation}
Hence the solutions $u_{0},u_{1}$ are actually related to the periods of the elliptic curve family
which are in turn related to modular forms for $\Gamma_{0}(2)$. 
This connection can then be used to give another proof of \eqref{eqnSchwarzianK3} by using modular forms.

Since our intention is to study the properties by using the Picard-Fuchs equations only, we shall not discuss the details. We wish to address the modularity and the application in Gromov-Witten theory of K3 orbifold surfaces somewhere else.
It is however worth pointing out that the symmetric square structure follows from the fact that the K3 surface family in consideration is polarized by some special lattices. Furthermore the K3 surfaces are Hodge-theoretically isomorphic to the Kummer varieties constructed from the elliptic curve family in \eqref{eqnE7family}. See \cite{Dolgachev:1996xw, Dolgachev:2013} for details.
\qed
\end{rem}

The symmetric square structure is essentially what leads to all these results,
it is therefore natural to 
ask when a third order ODE is the symmetric square of a second order one. This is an independent question on differential operators.

For a second order ODE
\begin{equation*}
\mathcal{L}_{2}= a_{2}\partial_{z}^{2}+a_{1}\partial_{z}+a_{0}\,,
\end{equation*}
its symmetric square is computed to be 
\begin{equation*}
\mathcal{L}_{3}= a_{2}^2\partial_{z}^{3}+3a_{1}a_{2}\partial_{z}^2+(a_{2}(a_{0}+\partial_{z}a_{1})+a_{1}
(2a_{1}-\partial_{z}a_{2})) \partial_{z}+(2a_{2}\partial_{z}a_{0}-2a_{0}\partial_{z}a_{2}+4a_{0}a_{1})  \,.
\end{equation*}
Hence a third order ODE which admits a symmetric square structure must have the form displayed above.
The condition of being a symmetric square can be more conveniently phrased in terms of the coefficients in the normal form of the differential operator, see \cite{Lian:1994zv} for detailed discussions.

\subsection{CY threefolds}

As discussed above, for K3 surfaces, the symmetric square structure implies \eqref{eqnSchwarzianK3} and leads to the triviality of genus zero Gromov-Witten invariants. 
One might wonder whether the same thing can be said for CY threefolds. \\

Consider for simplicity the one-parameter case.
We assume the family in the B-model $\check{\pi}:\check{\mathcal{X}}\rightarrow \check{\mathcal{B}}$
has a nice description in terms of algebraic varieties.
We also choose a local coordinate system $z$ on $\check{\mathcal{B}}$
in which the equation of the CY family and hence the Picard-Fuchs system is naturally written.

Under mirror symmetry, the "large volume limit" $t=i\infty$ in the space of K\"ahler structures, around which the Gromov-Witten theory is defined, is mapped to the so-called "large complex structure limit" \cite{Morrison:1993} which we assume is given by $z=0$. The mirror map then induces a new local coordinate system near the large complex structure limit which we also denote by $t$.\\

The base $\check{\mathcal{B}}$ is equipped with the Weil-Petersson metric whose K\"ahler potential $K$ is defined by
\begin{equation}\label{eqnKahlerpotential}
e^{-K(z,\bar{z})}=\sqrt{-1}\int_{\check{\mathcal{X}}_{z}} \Omega(z)\wedge \overline{\Omega(z)}\,,
\end{equation}
where $\Omega$ is a local holomorphic section of the Hodge line bundle over the base $\check{\mathcal{B}}$.
The Weil-Petersson geometry has many nice properties known as special geometry \cite{Strominger:1990pd}.
One property is that near\footnote{The structure in \eqref{eqnsymplecticperiods} for the periods actually holds everywhere on the base.} the large complex structure limit point $z$ on $\check{\mathcal{B}}$, there exists a holomorphic function $F(t)$ called prepotential so that 
the periods with respect to a symplectic basis of $H_{3}(\check{\mathcal{X}}_{z},\mathbb{Z})$ has the structure
\begin{equation}\label{eqnsymplecticperiods}
(\pi_{0}(z),\pi_{1}(z),\pi_{2}(z),\pi_{3}(z))=\pi_{0}(z)(1,t,\partial_{t}F,2F-t\,\partial_{t}F)\,.
\end{equation}
Moreover, the Yukawa coupling, as the mirror counterpart of the generating series of genus zero Gromov-Witten invariants of the CY family $\mathcal{X}$ in the A-model, is related to the prepotential $F(t)$ by
\begin{equation}
C_{ttt}=\partial_{t}^{3}F(t)\,.
\end{equation}
The prepotential has the following general form, see \cite{Cox:2000vi} and references therein,
\begin{equation}\label{eqnprepotential}
F(t)={\kappa\over 3!}t^{3}+c_{2}t^2+ c_{1}  t+ {\chi\over 2}\zeta(3)+F_{\mathrm{inst}}(q)
:={\kappa\over 3!}t^{3}+Q_{2}(t)+F_{\mathrm{inst}}(q)\,, \quad q=e^{-t}\,,
\end{equation}
where $\kappa, c_{1}, c_{2},\chi$ are some real constants depending on the family $\mathcal{X}$.
Now we can immediately see that the existence of a symmetric cubic structure on the Picard-Fuchs equation is equivalent to
the statement that $F(t)$ is cubic in $t$, or 
\begin{equation}
C_{ttt}=\kappa\,.
\end{equation}
That is, all of the Gromov-Witten invariants of non-zero degree vanish.
This is rarely the case for CY threefolds. See \cite{Ceresole:1992su, Ceresole:1993qq} for examples and further discussions.\\

Similar to the symmetric square case, the criteria for the symmetric cubic structure can be phrased in terms of some invariants constructed out of the coefficients of the ODE,
see \cite{Lian:1994zv, Ceresole:1992su, Ceresole:1993qq}.

\subsection{Quantum correction in Weil-Petersson geometry}

The Weil-Petersson metric for any CY family is defined through the same formula in \eqref{eqnKahlerpotential}.
It is easy to see that for elliptic curve families or the one-parameter K3 families admitting the symmetric square structure, the normalized period
takes values in the upper-half plane and the Weil-Petersson metric is exactly the Poincar\'e metric.\\

For CY threefold families, by restricting to an one-dimensional slice, one has
\begin{equation*}
e^{-K}={1\over 6}\kappa (t-\bar{t})^{3}
+(t-\bar{t})(\partial_{t}Q_{2}+\partial_{t}F_{\mathrm{inst}}+\overline{\partial_{t}Q_{2}+\partial_{t}F_{\mathrm{inst}}})-(2Q_{2}+2F_{\mathrm{inst}}-
\overline{2Q_{2}+2F_{\mathrm{inst}}})\,.
\end{equation*}
This gives rise to the Poincar\'e metric if and only if the summation of the terms on the right hand side of the above expression, except for the cubic term in $(t-\bar{t})$, is zero.
This implies $F_{\mathrm{inst}}=0$ or equivalently $C_{ttt}=\kappa$.
That is, for an one-parameter family, the Weil-Petersson metric on the base of the CY family is the quantum correction (by genus zero Gromov-Witten invariants)
of the Poincar\'e metric.
This is expected by Schmid's $\mathrm{SL}_{2}$-orbit theorem \cite{Schmid:1973variation}.
Similar statements also hold for multi-parameter families. 

Interestingly, from the Weil-Petersson metric on the base, one can define differential rings \cite{Yamaguchi:2004bt, Alim:2007qj, Hosono:2008ve, Alim:2013eja, Zhou:2013hpa} 
which exhibit similar structures as the rings of quasi-modular and almost-holomorphic modular forms \cite{Kaneko:1995} defined from the Poincar\'e metric and seem to provide generalizations thereof.

\newcommand{\etalchar}[1]{$^{#1}$}
\providecommand{\bysame}{\leavevmode\hbox to3em{\hrulefill}\thinspace}
\providecommand{\MR}{\relax\ifhmode\unskip\space\fi MR }
\providecommand{\MRhref}[2]{%
  \href{http://www.ams.org/mathscinet-getitem?mr=#1}{#2}
}
\providecommand{\href}[2]{#2}

\medskip{}
\noindent{\small Perimeter Institute for Theoretical Physics, 31 Caroline Street North, Waterloo, ON N2L 2Y5, Canada}
\noindent{\small E-mail: \tt jzhou@perimeterinstitute.ca}

\end{document}